\documentclass{elsart} 

\usepackage{amsmath,amssymb} 

\newtheorem{relation}[thm]{Relation}

\DeclareMathOperator{\Hom}{Hom}
\newcommand{\Q}{{\mathbb Q}}
\newcommand{\R}{{\mathbb R}}

\newcommand{\B}{{\mathfrak B}}
\newcommand{\bB}{{\mathbf B}}
\newcommand{\bF}{{\mathbf F}}

\begin{document}

\begin{frontmatter}

\title{Algebraic Structures of Bernoulli Numbers and Polynomials}

\author{I-Chiau Huang}
\address{Institute of Mathematics, Academia Sinica, Nankang, Taipei 11529, Taiwan, R.O.C.}
\ead{ichuang@math.sinica.edu.tw}

\begin{abstract}
In the field of Laurent series $\Q((T))$, we construct a subring $\B$ which has a natural 
$D$-module structure. Identities of Bernoulli numbers and polynomials are obtained from the algebraic structures of $\B$.
\end{abstract}


\end{frontmatter}


\section{Introduction}


Let $e^T=\sum_{i=0}^{\infty}T^i/i!$ be the formal power series with coefficients in the field
$\Q$ of rational numbers. In the field of Laurent series ${\mathbb Q}((T))$, the series
$\bB=\bB(T)=T/(e^T-1)$ is contained  in the formal power series ring ${\mathbb Q}[[T]]$. 
The $i$-th Bernoulli number $B^{(n)}_i$ of  order $n$ is defined by
$$
\bB^n=\sum_{i=0}^{\infty}\frac{B^{(n)}_i}{i!}T^i.
$$
We write $B^{(1)}_i$ simply as $B_i$. 
The $i$-th Bernoulli polynomial $B^{(n)}_i(X)$ of order $n$ is defined by
$$
\bB^ne^{XT}=
\sum_{i=0}^{\infty}\frac{B^{(n)}_i(X)}{i!}T^i
$$
in the power series ring $\Q[X][[T]]$, whose coefficients are in the polynomial ring 
$\Q[X]$.

Identities of Bernoulli numbers are abundant ever since early days of their discovery and still of contemporary interests. Among diverse approaches to these identities, we would like to explore 
their algebraic structures. These identities, often of order one, are better understood through 
Bernoulli numbers and polynomials of higher orders. A framework for such a phenomenon is built 
up by the vector spaces defined below.
\begin{defn}
We define $\B$ to be the $\Q$-vector subspace of $\Q((T))$ generated by Laurent series
$T^{m}\bB^n(bT)e^{aT}$, where $m$ is an integer, $n$ is a non-negative integer, 
$a\in\mathbb Q$ and $0\neq b\in\mathbb Q$. 
We define $\B_0$ to be the $\Q$-vector subspace of $\B$ generated by $T^{m}\bB^n$. 
\end{defn}

The power series $\bB^n(bT)e^{aT}\in\B$ is the exponential generating function of 
$b^iB^{(n)}_i(a/b)$, which is the main ingredient of our identities. A multiplication on 
$\bB^n(bT)e^{aT}$ by a power of $T$ simply means a shifting of the generating function. 
These observations motivate our definition.

Clearly $\B_0$ is a subring of $\Q((T))$. In Proposition \ref{prop:ring}, we will show 
that $\B$ is a subring as well. Note that $\bB^{-1}=T^{-1}e^T-T^{-1}\in\B$. Hence, for 
a fixed $k$, the exponential generating function of the sequence 
$\{S(n,k)\}_{n\geq 0}$ of Stirling numbers of the second kind is contained in $\B$. Indeed
$$
\sum_{n=0}^{\infty}S(n,k)\frac{T^n}{n!}=\frac{T^k}{k!}\bB^{-k}.
$$
Besides the $\Q[T,e^T]$-module structure, the field $\Q((T))$ has a natural left module 
structure over the Weyl algebra $D:=\Q\langle T,d/dT\rangle$. It is clear that $\B$ is 
a $\Q[T,e^T]$-submodule of $\Q((T))$. In Proposition \ref{prop:Dmodule}, we will show 
that $\B$ is also a $D$-submodule of $\Q((T))$. The proof is essentially due to 
Lucas \cite[p. 62]{luc:ds}. In terms of our terminology, he showed that $\B_0$ is a 
$D$-submodule of $\Q((T))$.

Starting with elements of $\B$ and performing algebraic operations on them, we will 
provide methods to represent the result as a $D$-linear combination of elements of the form 
$T^{m}\bB(bT)e^{aT}$. Equating coefficients in the Laurent series of the representation, we 
obtain identities of Bernoulli numbers and Bernoulli polynomials. For instance, the relation
$$
\bB^2=\left(1-T-T\frac{d}{dT}\right)\bB
$$
gives rise to Euler's identity
\begin{equation}\label{eq:Euler}
\sum_{i=1}^{m-1}\binom{2m}{2i}B_{2i}B_{2m-2i}=-(2m+1)B_{2m}
\end{equation}
by equating the coefficients of $T^{2m}$ for $m>1$. See Relation~\ref{iden:lo} for details. 
Such a viewpoint applies to many identities found in the literature including generalizations 
of Euler's identity \cite{dil:spbn,san:iirzf,sit-dav:siirzf,zhan:sirzf}, an identity from 
Eisenstein series \cite{rad:tant}, identities from special values of Zeta 
functions \cite{eie:nbnsgbp,eie-lai:bia}, and others \cite{ago-dil:cilrbn}.
One may implement algebraic operations of $\B$ into computer programs, which provide new 
identities elegant or involved.

It is plausible that our algebraic view may be extended to other classical numbers. 
We leave the possibilities to the reader.


\section{Algebraic Structures}


Note that elements $T^{m}\bB^n$ form a basis for the $\Q$-vector space $\B_0$.
This fact is a consequence of the transcendence of $\bB$ over $\Q[T]$, which is 
equivalent to the transcendence of $e^T$. In the $\Q[T,e^T]$-module $\B$, elements 
$T^{m}\bB^n(bT)e^{aT}$ may have non-trivial relations. For instance,
$$
\bB e^T=T+\bB.
$$
See Section \ref{sec:iden} for more examples.

The ring structure of $\B$ is based on the following version of partial fraction decompositions.
\begin{lem}\label{lem:941145}
For any positive integers $n_1,\cdots,n_s,k_1,\cdots,k_s$, there exist
$g_i\in{\mathbb Q}[X]$ and positive integers
$m_i$, $\ell_i$ with $\ell_i\leq n_1+\cdots+n_s$ such
that
\begin{equation}\label{eq:59911}
\frac{1}{(X^{k_1}-1)^{n_1}\cdots(X^{k_s}-1)^{n_s}}=
\sum_{i}\frac{g_i}{(X^{m_i}-1)^{\ell_i}}.
\end{equation}
\end{lem}
\begin{pf}
We may assume that $k_1,\cdots,k_s$ are distinct.
The case $s=1$ is trivial. Assume that $s=2$. Let $k$ be the greatest common divisor of $k_1$
and $k_2$. Since $X-1$, $X^{(k_1/k)-1}+\cdots+X+1$ and $X^{(k_2/k)-1}+\cdots+X+1$
are relatively prime, there exist $g^*_0,g^*_1,g^*_2\in{\mathbb Q}[X]$ such that
\begin{small}\begin{eqnarray*}
& & \frac{1}{(X^{k_1/k}-1)^{n_1}(X^{k_2/k}-1)^{n_2}}\\
&=&
\frac{g^*_0}{(X-1)^{n_1+n_2}}+
\frac{g^*_1}{(X^{(k_1/k)-1}+\cdots+X+1)^{n_1}}+
\frac{g^*_2}{(X^{(k_2/k)-1}+\cdots+X+1)^{n_2}}.
\end{eqnarray*}\end{small}
Let $g_0=g^*_0(X^k)$, $g_1=(X^k-1)^{n_1}g^*_1(X^k)$,
$g_2=(X^k-1)^{n_2}g^*_2(X^k)$. Then
$$
\frac{1}{(X^{k_1}-1)^{n_1}(X^{k_2}-1)^{n_2}}=
\frac{g_0}{(X^k-1)^{n_1+n_2}}+
\frac{g_1}{(X^{k_1}-1)^{n_1}}+
\frac{g_2}{(X^{k_2}-1)^{n_2}}.
$$
The lemma follows from induction on $s$.
\end{pf}

Since the denominators of (\ref{eq:59911}) has a common factor $X-1$, the decomposition in 
Lemma~\ref{lem:941145} is not unique. For instance,
$$
\frac{1}{(X-1)(X^2-1)}=
\frac{1-\frac{1}{2}X}{(X-1)^2}+\frac{\frac{1}{2}X}{X^2-1}=
\frac{\frac{1}{2}}{(X-1)^2}+\frac{-\frac{1}{2}}{X^2-1}.
$$

\begin{prop}\label{prop:ring}
$\B$ is a subring of $\Q((T))$.
\end{prop}
\begin{pf}
We need to show that 
$\bB^{n_1}((k_1/k)T)\bB^{n_2}((k_2/k)T)\in\B$, where $k_1,k_2$ are positive integers ad $k$ is an
integer. Choose polynomials $g_i$ and positive  integers $m_i$, $\ell_i$ as in 
Lemma~\ref{lem:941145} for $s=2$. Replacing
$X$ by $e^{T/k}$ in (\ref{eq:59911}) and multiplying both sides of the equation by 
$(k_1T)^{n_1}(k_2T)^{n_2}/k^{n_1+n_2}$, we obtain the required formula
$$
\bB^{n_1}(\frac{k_1}{k}T)\bB^{n_2}(\frac{k_2}{k}T)=
\sum_{i}(\frac{k_1}{k})^{n_1}(\frac{k_2}{k})^{n_2}(\frac{k}{m_i})^{\ell_i}
T^{n_1+n_2-\ell_i}\bB^{\ell_i}(\frac{m_i}{k}T)g_i(e^{T/k}).
$$
\end{pf}

The Weyl algebra $D=\Q\langle T,d/dT\rangle$ is generated by the multiplication operator $T$ and the derivation $d/dT$ as a $\Q$-subalgebra of $\Hom_\Q(\Q((T)),\Q((T)))$ . For $f,g\in D$, we write 
their product as $f\cdot g$. The field of Laurent 
series $\Q((T))$ is endowed with the natural left $D$-module structure. For instance, 
the module structure gives $(d/dT)T=1\in\Q((T))$. Do not be confused with the product 
$(d/dT)\cdot T$, which satisfies the Leibniz rule
$$
\frac{d}{dT}\cdot T=1+T\cdot\frac{d}{dT}.
$$
The subalgebra of $D$ generated by $T$ is also denoted by $\Q[T]$, which makes no confusion with
the subring of $\Q((T))$ generated by $T$.
For $f,g\in\Q[T]\subset D$, their product is also written as $fg$. This notation agrees with 
that of the module structure on $\Q((T))$. Moreover, using the Leibniz rule,
any element of $D$ can be written as $\sum f_i\cdot(d/dT)^i$, where $f_i\in\Q[T]$. 
We write $\sum f_i\cdot(d/dT)^i$ also as $\sum f_i(d/dT)^i$
without any confusion with the module structure.
More notations: $d^i/dT^i$ stands for $(d/dT)^i$.  For $\bF\in\Q((T))$, we use 
$d\bF^n/dT^n$ for $(d/dT)^n\bF$. In our convention,
$d\bF^0/dT^0=\bF$.

\begin{prop}\label{prop:Dmodule}
$\B_0$ and $\B$ are $D$-submodules of $\Q((T))$. As $D$-modules, $\B$ is generated by 
$T^m\bB(bT)e^{aT}$ and $\B_0$ is generated by $T^m\bB$.
\end{prop}
\begin{pf}
For the first assertion, it suffices to show $(d/dT)(\bB^n(bT)e^{aT})\in\B$ and
$(d/dT)\bB^n\in\B_0$.  This follows from the formula
\begin{equation}\label{eq:Lucas}
\frac{d}{dT}\left(\bB^n(bT)e^{aT}\right)=
\left((a+\frac{n}{T})\bB^n(bT)-nb\bB^n(bT)-
\frac{n}{T}\bB^{n+1}(bT)\right)e^{aT},
\end{equation}
which can be found in \cite[p. 62]{luc:ds} for $b=1$ and $a=0$. The general case of
(\ref{eq:Lucas}) reduces easily to the case $a=0$, which is obtained by dividing
the equality
\begin{eqnarray*}
& & n\bB^n(bT)(e^{bT}-1)^n \\ 
&=& n(bT)^n \\
&=& T\frac{d(bT)^n}{dT}\\
&=& T\frac{d}{dT}\left(\bB^n(bT)(e^{bT}-1)^n\right)\\
&=& nbT\bB^n(bT)(e^{bT}-1)^{n-1}e^{bT}+
        T(e^{bT}-1)^n\frac{d}{dT}\bB^n(bT) \\
&=& nbT\bB^n(bT)(e^{bT}-1)^n+
        n\bB^{n+1}(bT)(e^{bT}-1)^n+T(e^{bT}-1)^n\frac{d}{dT}\bB^n(bT)
\end{eqnarray*}
by $T(e^{bT}-1)^n$, {\em c.f.} \cite[Prop. 1]{hu-hu:bnpr}.
The second assertion follows from
\begin{equation}\label{eq:877199}
\bB^{n+1}(bT)e^{aT}=\left(1-bT+\frac{aT}{n}-\frac{T}{n}\frac{d}{dT}\right)
\left(\bB^n(bT)e^{aT}\right),
\end{equation}
which is equivalent to (\ref{eq:Lucas}).
\end{pf}

In Section~\ref{sec:iden}, we will provide two special cases of partial fraction 
decompositions which suffice to calculate the product of elements in $\B$. 

The first case is a decomposition of $1/(X^n-1)(X^m-1)$, where $m\neq n$.
Let $\ell$ be the greatest common divisor of $m$ and $n$. 
\begin{defn}\label{defn:g_mn}
For $\ell=1$, we define $g_{m,n}$ to be the polynomial with degree less than $m-1$ such that
\begin{eqnarray*}
&&
\frac{1}{(X-1)(1+X+\cdots+X^{m-1})(1+X+\cdots+X^{n-1})}\\
&=&\frac{1}{mn(X-1)}+\frac{g_{n,m}}{1+X+\cdots+X^{n-1}}+
\frac{g_{m,n}}{1+X+\cdots+X^{m-1}}.
\end{eqnarray*}
In general, we define $g_{n,m}=g_{n/\ell,m/\ell}(X^\ell)$.
\end{defn}
For example,
\begin{eqnarray*}
g_{2,3}=-\frac{1}{2}, && g_{3,2}=\frac{1}{3}(X-1);\\
g_{2,5}=-\frac{1}{2}, && g_{5,2}=\frac{1}{5}(2X^3-X^2+X-2);\\
g_{3,5}=\frac{1}{3}(X-1), && g_{5,3}=\frac{1}{5}(-2X^3+X^2-X-3).
\end{eqnarray*}
We have a decomposition
$$
\frac{1}{(X^n-1)(X^m-1)}=
\frac{\ell^2}{mn(X^\ell-1)^2}+\frac{g_{n,m}}{X^n-1}+\frac{g_{m,n}}{X^m-1}.
$$
Note that $g_{1,m}=0$. For $n>1$ and $\ell=1$, we may describe
$g_{n,m}=a_0+a_1X+\cdots+a_{n-2}X^{n-2}$ 
using the complex number $v=e^{2\pi\sqrt{-1}/n}$ by
$$
\begin{pmatrix}a_0\\a_1\\ \vdots\\a_{n-2}\end{pmatrix}
=
\begin{pmatrix}1&v&\cdots&v^{n-2}\\1&v^2&\cdots&v^{2n-4}\\
\vdots&\vdots&\ddots&\vdots\\
1&v^{n-1}&\cdots&v^{(n-1)(n-2)}\end{pmatrix}^{-1}
\begin{pmatrix}(v^m-1)^{-1}\\(v^{2m}-1)^{-1}\\ 
\vdots\\(v^{(n-1)m}-1)^{-1}\end{pmatrix}.
$$

The second case is a decomposition of $1/(X^\ell-1)^k(X^n-1)$, where $\ell$ is a divisor of $n$.
\begin{defn}\label{defn:fh}
We define $h^{(k)}_{1,n}$ and $f^{(k)}_{1,n}$ to be the polynomials with degrees less than
$k$ and $n-1$ respectively such that
$$
\frac{1}{(X-1)^k(1+X+\cdots+X^{n-1})}=
\frac{h^{(k)}_{1,n}}{(X-1)^k}+\frac{f^{(k)}_{1,n}}{1+X+\cdots+X^{n-1}}.
$$
We define $h^{(k)}_{\ell,n}=h^{(k)}_{1,n/\ell}(X^\ell)$ and 
$f^{(k)}_{\ell,n}=f^{(k)}_{1,n/\ell}(X^\ell)$.
\end{defn}
For example,
\begin{eqnarray*}
h^{(2)}_{1,5}&=&\frac{1}{5}(-2X+3),\\
f^{(2)}_{1,5}&=&\frac{1}{5}(2X^3+3X^2+3X+2).
\end{eqnarray*}
We have a decomposition
$$
\frac{1}{(X^\ell-1)^k(X^n-1)}=
\frac{h^{(k)}_{\ell,n}}{(X^\ell-1)^{k+1}}+\frac{f^{(k)}_{\ell,n}}{X^n-1}.
$$
Note that $h^{(k)}_{1,n}=a_0+a_1(X-1)+\cdots+a_{k-1}(X-1)^{k-1}$ can be described inductively 
by $a_0=1/n$ and
$$
a_0\binom{n}{i}+a_1\binom{n}{i-1}+\cdots+a_{i-1}\binom{n}{1}=0
$$
for $i>1$ with our convention $\binom{n}{i}=0$ for $n<i$. Note also that
$$
f^{(k)}_{1,n}=\frac{1-(1+X+\cdots+X^{n-1})h^{(k)}_{1,n}}{(X-1)^k}.
$$


\section{Relations and Identities}\label{sec:iden}


We list some relations in $\B$ and demonstrate their implications on identities 
of Bernoulli numbers and polynomials.


\begin{relation}\label{rel:1}
$$
\bB(-T)-T=\bB.
$$
\end{relation}
It is straightforward to verify the above relation, which implies that $B_1=-1/2$ 
and all other odd Bernoulli numbers of order one are zero.


\begin{relation}\label{rel:2}
$$
e^{T}\bB=\bB(-T).
$$
\end{relation}
It is straightforward to verify the above relation, which implies the recurrence formula 
$$
\sum_{i=0}^n\binom{n}{i}B_i=(-1)^nB_n
$$
for $n\geq 0$. Together with the vanishing property of odd Bernoulli numbers of order one, the above recurrence formula can be also written as
$$
\sum_{i=0}^n\binom{n}{i}B_i=B_n
$$
for $n\geq 2$. See \cite[p.18, Equation (4)]{nor:d}. 


\begin{relation}\label{iden:mul}
$$
\sum_{i=0}^{n-1}\bB(nT)e^{iT}=n\bB.
$$
\end{relation}
Relation~\ref{iden:mul} follows from the direct computation
$$
\sum_{i=0}^{n-1}\bB(nT)e^{iT}
=n\frac{T}{e^T-1}\frac{e^{nT}-1}{nT}\bB(nT)
=n\bB.
$$
It is equivalent to
$\sum_{i=0}^{n-1}\bB e^{(a+\frac{i}{n})T}=n\bB(\frac{T}{n})e^{aT}$, whose 
coefficients of $T^m$ give rise to multiplication theorem
$$
\sum_{i=0}^{n-1}B_m(a+\frac{i}{n})=n^{1-m}B_m(na)
$$
\cite[p. 21, Equation (18)]{nor:d}.


\begin{relation}\label{iden:lo}
$$
\bB^{n+1}(bT)e^{aT}=
\left(1-bT+\frac{aT}{n}-\frac{T}{n}\frac{d}{dT}\right)\left(\bB^n(bT)e^{aT}\right).
$$
\end{relation}

This is equation~(\ref{eq:877199}).
For $b=1$, equating its coefficients of $T^i$, we get the formula for lowering orders of 
Bernoulli polynomials
$$
B^{(n+1)}_i(a)=(1-\frac{i}{n})B^{(n)}_i(a)+(a-n)\frac{i}{n}B^{(n)}_{i-1}(a)
$$
\cite[p. 145, Equation (81)]{nor:d}.

For $a=0$ and $b=n=1$, we get Euler's identity
$$
\sum_{i=0}^{m}\binom{m}{i}B_iB_{m-i}=(1-m)B_m-mB_{m-1}
$$
from the coefficients of $T^m$.
Using the vanishing property for odd Bernoulli numbers of order one, we recast Euler's
identity in the form of (\ref{eq:Euler}).

For $a=0$, $b=1$ and $n>1$, repeatedly using Relation~\ref{iden:lo}, we get generalizations of 
Euler's identity. See R.~Sitaramachandrarao and B.~Davis \cite{sit-dav:siirzf}, 
A.~Sankarayanan \cite{san:iirzf},  W.-P.~Zhang \cite{zhan:sirzf}, and 
K.~Dilcher \cite{dil:spbn}. 

The polynomial analogue \cite[Proposition 2]{eie:nbnsgbp}
$$
\sum_{i=0}^n\binom{n}{i}B_i(a)B_{n-i}(b)=(1-n)B_n(a+b)+n(a+b-1)B_{n-1}(a+b)
$$
of (\ref{eq:Euler}) also follows from direct 
computation:
\begin{eqnarray*}
(\bB e^{aT})(\bB e^{bT})&=&\left((1-T-T\frac{d}{dT})\bB\right)e^{(a+b)T}\\
&=&\left(1+(a+b-1)T-T\frac{d}{dT}\right)\left(\bB e^{(a+b)T}\right).
\end{eqnarray*}


\begin{relation}\label{iden:cs}
$$
\frac{d^n\bB}{dT^n}=T^{-n}f_n(T,\bB),
$$
where 
$$
f_n:=(-1)^n\sum_{j=1}^{n+1}(j-1)!\left(S(n+1,j)U^{n-j+1}-nS(n,j)U^{n-j}\right)V^j
$$
with Stirling numbers $S(n,j)$ of the second kind.
\end{relation}
Relation~\ref{iden:cs} is given in \cite[Theorem 1]{hu-hu:bnpr} with
$f_n\in{\mathbb Z}[U,V]$ constructed inductively by $f_0:=V$ and 
$$
f_n:=\left(1-n+U\frac{\partial}{\partial U}+V\frac{\partial}{\partial V}
-UV\frac{\partial}{\partial V}-V^2\frac{\partial}{\partial V}\right)f_{n-1}
$$ 
for $n>0$. Using the initial values $S(n,1)=S(n,n)=1$, $S(n,n+1)=0$ and the recurrence relation 
$$
S(n+1,j)=S(n,j-1)+jS(n,j),
$$
one shows that the two definitions for $f_n$ agree, {\em cf.} 
\cite[Lemma 3.1]{ago-dil:cilrbn}. 

Given $m$ and $n$, 
$$
T^{m+n}\frac{d^m\bB}{dT^m}\frac{d^n\bB}{dT^n}=f_m(T,\bB)f_n(T,\bB)
=\varphi\bB
$$ 
for some $\varphi\in D$ by Relation~\ref{iden:lo}. In terms of Bernoulli numbers, this 
gives rise to identities described by T. Agoh and 
K. Dilcher \cite{ago-dil:cilrbn}. For example, the identity
$$
\sum_{i=0}^n\binom{n}{i}B_{1+i}B_{1+n-i}=
\frac{1}{6}(n-1)B_{n}-B_{n+1}-\frac{1}{6}(n+3)B_{n+2}
$$
comes from the relation
$$
\frac{d\bB}{dT}\frac{d\bB}{dT}
=
\left(-\frac{1}{6}T\frac{d^3}{dT^3}-\frac{1}{2}\frac{d^2}{dT^2}+
\frac{1}{6}T\frac{d}{dT}-\frac{d}{dT}-\frac{1}{6}\right)\bB.
$$

Combining with Relation \ref{iden:lo}, we have another relation
$$
T^2\left(\frac{d\bB}{dT}\right)^2-(2n-1)\bB^2=\varphi\bB,
$$
where 
\begin{eqnarray*}
\varphi&=&\left(-\frac{1}{6}T^3\frac{d^3}{dT^3}-\frac{1}{2}T^2\frac{d^2}{dT^2}+
\frac{1}{6}T^3\frac{d}{dT}-T^2\frac{d}{dT}-\frac{1}{6}T^2\right)\\
&&
-(2n-1)\left(1-T-T\frac{d}{dT}\right).
\end{eqnarray*}
Equating coefficients of $T^{2n}$ with $n\ge 4$, we get the identity
$$
\sum_{i=2}^{n-2}
\frac{(2n-2)!}{(2i-2)!(2n-2i-2)!}\frac{B_{2i}}{2i}\frac{B_{2n-2i}}{2n-2i}
=-\frac{(2n+1)(n-3)}{6n}B_{2n},
$$
which was proved by H.~Rademacher \cite{rad:tant} using Eisenstein series
and by M.~Eie \cite{eie:nbnsgbp} using Zeta functions. See also \cite{hu-hu:bnpr}.


\begin{relation}\label{iden:prod_mn}
$$
\bB(mT)\bB(nT)=\bB^2(\ell T)+mT\bB(nT)g_{n,m}(e^T)+nT\bB(mT)g_{m,n}(e^T),
$$
where 
$\ell$ is the  greatest common divisor of distinct positive integers $m,n$ and polynomials 
$g_{m,n}$, $g_{n,m}$ are given in Definition~\ref{defn:g_mn}.
\end{relation}
For example,
\begin{eqnarray*}
\bB(2T)\bB(3T)&=&\bB^2+\frac{2}{3}T(e^T-1)\bB(3T)-\frac{3}{2}T\bB(2T),\\
\bB(2T)\bB(5T)&=&\bB^2+\frac{2}{5}T(2e^{3T}-e^{2T}+e^T-2)\bB(5T)-\frac{5}{2}T\bB(2T),\\
\bB(3T)\bB(5T)&=&\bB^2+\frac{3}{5}T(-2e^{3T}+e^{2T}-e^T-3)\bB(5T)+\frac{5}{3}T(e^T-1)\bB(3T).
\end{eqnarray*}
In terms of Bernoulli polynomials, the above expression for $\bB(2T)\bB(3T)$ together with
the relation $\bB^2=(1-T-T(d/dT))\bB$ gives rise to
$$
\sum_{i=0}^n3^i2^{n-i}\binom{n}{i}B_iB_{n-i}
=
\frac{2n}{3^{2-n}}B_{n-1}(\frac{1}{3})-
\left(\frac{2n}{3^{2-n}}+\frac{3n}{2^{2-n}}+n\right)B_{n-1}+(1-n)B_n.
$$
For even Bernoulli numbers, the identity reduces to
$$
\sum_{i=0}^{n}3^{2i}2^{2n-2i}\binom{2n}{2i}B_{2i}B_{2n-2i}
=
\frac{4n}{3^{2-2n}}B_{2n-1}(\frac{1}{3})+(1-2n)B_{2n},
$$
where $n>1$.


\begin{relation}\label{iden:prod_nkl}
$$
\bB^k(\ell T)\bB(nT)=\ell^{k}T^{k}\bB(nT)f^{(k)}_{\ell,n}(e^T)+
\dfrac{n}{\ell}\bB^{k+1}(\ell T)h^{(k)}_{\ell,n}(e^T),
$$ 
where $\ell$ is a divisor of $n$ and $f^{(k)}_{\ell,n}$, $h^{(k)}_{\ell,n}$ are given in
Definition~\ref{defn:fh}.
\end{relation}
For example,
$$
\bB^2\bB(5T)=\frac{1}{5}T^2(2e^{3T}+3e^{2T}+3e^T+2)\bB(5T)+(-2e^T+3)\bB^3.
$$
Relation~\ref{iden:prod_nkl} writes $\bB^k(\ell T)\bB(nT)$ as a
$\Q[T,e^T]$-linear combination of $\bB(nT)$ and $\bB^{k+1}(\ell T)$.
Multiply Relation~\ref{iden:prod_mn} by $\bB(nT)$ and using Relation~\ref{iden:prod_nkl}, 
we are 
able to write $\bB(mT)\bB^2(nT)$ as a $\Q[T,e^T]$-linear combination of $\bB(nT)$, 
$\bB^2(nT)$, $\bB^2(\ell T)$, $\bB^3(\ell T)$ and $\bB(mT)$. Repeating the process, we are 
able to write $\bB(mT)\bB^k(nT)$ as a $\Q[T,e^T]$-linear combination of $\bB(nT)$, 
$\cdots$, $\bB^k(nT)$, $\bB^2(\ell T)$, $\cdots$, $\bB^{k+1}(\ell T)$ and $\bB(mT)$. Thus 
we are able to write the product $\bB(m_1T)\cdots\bB(m_sT)$ as a $\Q[T]$-linear combination of
elements of the form $\bB^n(kT)e^{aT}$. Using Relation~\ref{iden:lo}, the product
$\bB(m_1T)\cdots\bB(m_sT)$ can be written as a $D$-linear combination of elements of the form
$\bB(bT)e^{aT}$.

Here we carry out the computation for $\bB(2T)\bB(3T)\bB(5T)$.
Relations \ref{iden:prod_mn} and \ref{iden:prod_nkl} together with the relation 
$(e^T-1)\bB^3=T\bB^2$ (from Relations~\ref{rel:1} and \ref{rel:2}) give rise to
\begin{eqnarray*}
&&
\bB(2T)\bB(3T)\bB(5T)\\
&=&\bB^2\bB(5T)+\frac{2}{3}T(e^T-1)\bB(3T)\bB(5T)-\frac{3}{2}T\bB(2T)\bB(5T)\\
&=&
\frac{1}{5}T^2(2e^{3T}+3e^{2T}+3e^T+2)\bB(5T)-2T\bB^2+\bB^3\\
&&
+\frac{2}{3}T(e^T-1)\bB^2+\frac{2}{5}T^2(e^T-1)(-2e^{3T}+e^{2T}-e^T-3)\bB(5T)\\
&&+\frac{10}{9}T^2(e^{2T}-2e^T+1)\bB(3T)
\\
&&
-\frac{3}{2}T\bB^2-\frac{3}{5}T^2(2e^{3T}-e^{2T}+e^T-2)\bB(5T)+\frac{15}{4}T^2\bB(2T).
\end{eqnarray*}
By Relation \ref{iden:mul},
\begin{eqnarray*}
(e^T-1)(-2e^{3T}+e^{2T}-e^T-3)\bB(5T)&=&5(e^{3T}+1)\bB(5T)-10\bB,\\
(e^{2T}-2e^T+1)\bB(3T)&=&-3e^T\bB(3T)+3\bB.
\end{eqnarray*}
Together with the relation $(e^T-1)\bB^2=T\bB$, we obtain
\begin{eqnarray*}
\bB(2T)\bB(3T)\bB(5T)
&=&
\bB^3-\frac{7}{2}T\bB^2+\frac{15}{4}T^2\bB(2T)-\frac{10}{3}T^2e^T\bB(3T)\\
&&
+\frac{6}{5}T^2(e^{3T}+e^{2T}+3)\bB(5T).
\end{eqnarray*}
By Relation \ref{iden:lo},
\begin{eqnarray*}
\bB(2T)\bB(3T)\bB(5T)
&=&
\left(\frac{1}{2}T^2\frac{d^2}{dT^2}+5T^2\frac{d}{dT}-T\frac{d}{dT}+
\frac{9}{2}T^2-5T+1\right)\bB\\
&&
+\frac{15}{4}T^2\bB(2T)-\frac{10}{3}T^2\bB(3T)e^T\\
&&
+\frac{6}{5}T^2\bB(5T)e^{3T}+\frac{6}{5}T^2\bB(5T)e^{2T}+\frac{18}{5}T^2\bB(5T).
\end{eqnarray*}
In terms of Bernoulli polynomials,
\begin{eqnarray*}
&&
\sum_{\stackrel{\scriptstyle i,j,k\geq 0}{i+j+k=n}}
\binom{n}{i,j,k}2^i3^j5^kB_iB_jB_k\\
&=&
\frac{1}{2}(n-1)(n-2)B_n+5n(n-2)B_{n-1}\\
&&
+\left(\frac{9}{2}
+\frac{15}{4} 2^{n-2}
+\frac{18}{5} 5^{n-2}\right)n(n-1)B_{n-2}-\frac{10}{3}n(n-1)3^{n-2}B_{n-2}(\frac{1}{3})\\
&&
+\frac{6}{5}n(n-1)5^{n-2}B_{n-2}(\frac{2}{5})+\frac{6}{5}n(n-1)5^{n-2}B_{n-2}(\frac{3}{5}).
\end{eqnarray*} 

Related to the representation of the product $\bB(2T)\bB(3T)\bB(5T)\bB(6T)$, an identity of
Bernoulli polynomials can be found in \cite[Proposition 4]{eie-lai:bia}.


\section{Variations}


There are identities of Bernoulli numbers, which do not come directly
from the algebraic structures of $\B$. For an example, we replace the underlying field
$\Q$ by the field $\R$ of real numbers.

\begin{relation}
For $s\in\R$,
$$
\bB(sT)\bB((1-s)T)=
(1-s)\left(\bB(sT)+\frac{sT}{2}\right)\bB+s\left(\bB((1-s)T)+\frac{(1-s)T}{2}\right)\bB.
$$
\end{relation}
This relation follows from direct computation and is equivalent to
\begin{eqnarray*}
&&
(\bB(sT)-1)(\bB((1-s)T)-1)\\
&=&
(1-s)\left(\bB(sT)+\frac{sT}{2}-1\right)\bB+
s\left(\bB((1-s)T)+\frac{(1-s)T}{2}-1\right)\bB\\
&&
+1+\bB-\bB(sT)-\bB((1-s)T).
\end{eqnarray*}
Equating the coefficients of $T^n$ for $n>0$, we obtain
\begin{eqnarray*}
\sum_{\stackrel{\scriptstyle i,j>0}{i+j=n}}
s^i(1-s)^j\frac{B_i}{i!}\frac{B_j}{j!}
&=&
\sum_{\stackrel{\scriptstyle k>0}{\ell+2k=n}}
\left((1-s)s^{2k}+s(1-s)^{2k}\right)
\frac{B_\ell}{\ell!}\frac{B_{2k}}{(2k)!}\\
&&
+\left(1-s^n-(1-s)^n\right)\frac{B_n}{n!}.
\end{eqnarray*}
Dividing this identity by $s(1-s)$, the integrations
\begin{eqnarray*}
\int_0^1 s^i(1-s)^j\frac{ds}{s(1-s)}&=&\frac{(i-1)!(j-1)!}{(i+j-1)!},\\
\int_0^1 (1-s^n-(1-s)^n)\frac{ds}{s(1-s)}&=&2\sum_{1\leq\ell<n}\frac{1}{\ell}=2H_{n-1}
\end{eqnarray*}
give rise to 

$$
\sum_{\stackrel{\scriptstyle i,j>0}{i+j=n}}
\frac{1}{(n-1)!}\frac{B_i}{i}\frac{B_j}{j}
=
\sum_{\stackrel{\scriptstyle k>0}{\ell+2k=n}}
\frac{1}{k}\frac{B_\ell}{\ell!}\frac{B_{2k}}{(2k)!}
+2H_{n-1}\frac{B_n}{n!}.
$$
For $n\geq 4$, multiplying the above equation by $(n-1)!$ and rearranging the summations, we obtain 
Miki's identity \cite{mi:rbbn}
\begin{eqnarray*}
\sum_{i=2}^{n-2}\frac{B_i}{i}\frac{B_{n-i}}{n-i}
&=&
\frac{2}{n}H_nB_n+
\frac{1}{n}\sum_{k=2}^{n-2}\binom{n}{k}
\left(B_{n-k}\frac{B_k}{k}+B_k\frac{B_{n-k}}{n-k}\right)\\
&=&
\frac{2}{n}H_nB_n+\sum_{k=2}^{n-2}\binom{n}{k}\frac{B_k}{k}\frac{B_{n-k}}{n-k}.
\end{eqnarray*}
The above proof of Miki's identity is taken from M.~C. Crabb \cite{cra:mgbni}, which is essentially a distillation 
of an argument given by G.~V. Dunne and C. Schubert  \cite{dun-schub:bniqft}. The argument works
also for the generalization given by I.~M. Gessel \cite{ges:mibn}.

Identities of Bernoulli numbers may come from coefficients in non-trivial elements in $\B$.
Consider
$$
\varphi=(2k+1)T^k\frac{d^k}{dT^k}+T^{k+1}\frac{d^{k+1}}{dT^{k+1}}.
$$
We claim that the coefficient of $T^{k+1}$ in
$e^T\varphi\bB$ is always zero. For a proof, we use the formula 
$$
\frac{d}{dT}\cdot T^k=kT^{k-1}+T^k\frac{d}{dT}
$$
to write $\varphi$ as the sum of $\varphi_1$ and
$\varphi_2$, where $\varphi_1=(k+1)T^k(d^k/dT^k)$ and $\varphi_2=T(d/dT)\cdot
T^k(d^k/dT^k)$. Note that the coefficients of $T^{k+1}$ in $(k+1)!e^T\varphi_1\bB$ and
$(k+1)!e^T\varphi_2\bB$ are $\sum\binom{k+1}{j}(k+1)B_{k+j}$ and $\sum\binom{k+1}{j}jB_{k+j}$, respectively. Thus the identity 
$$
\sum_{j=0}^{k+1}\binom{k+1}{j}(k+j+1)B_{k+j}=0
$$
discovered by M. Kaneko \cite{kan:rfbn} is equivalent to our claim.

\end{document}